\documentclass[12pt]{amsart}
\usepackage{amsmath, amssymb, amsthm, mathrsfs, amscd}

\newcommand{\col}{\colon \negthinspace}

\newcommand{\m}{{\mathfrak{m}}}

\DeclareMathOperator{\tor}{Tor}

\DeclareMathOperator{\ord}{ord}

\newtheorem{theorem}{Theorem}[section]
\newtheorem{cor}[theorem]{Corollary}
\newtheorem{lemma}[theorem]{Lemma}

\newtheorem*{nonumthm}{Theorem}
\newtheorem*{T1}{Theorem 1}

\theoremstyle{definition}
\newtheorem{definition}[theorem]{Definition}

\theoremstyle{remark}
\newtheorem{remark}[theorem]{Remark}

\newtheorem*{rmks}{Remarks}

\numberwithin{equation}{section}

\begin{document}

\title{Hilbert-Kunz Functions for Normal Rings}

\author{Craig Huneke}
\address{Department of Mathematics,
University of Kansas,
Lawrence, KS 66045}
\email{huneke@math.ukans.edu}
\urladdr{http://www.math.ukans.edu/\textasciitilde huneke}

\author{Moira A. McDermott}
\address{Department of Mathematics and Computer Science,
Gustavus Adolphus College,
Saint Peter, MN 56082}
\email{mmcdermo@gac.edu}
\urladdr{http://www.gustavus.edu/\textasciitilde mmcdermo}

\author{Paul Monsky}
\address{Department of Mathematics,
Brandeis University,
Waltham, MA 02254-9110}
\email{monsky@brandeis.edu}
\urladdr{http://www.math.brandeis.edu/homepages/monsky}

\date{\today}
\thanks{The first author was partially supported by NSF grant
DMS-0244405, and the second author was partially
supported by the NSF-AWM Mentoring Travel Grant.  We thank them for
their support.}

\bibliographystyle{alpha}

\begin{abstract}
Let $(R,\m,k)$ be an excellent, local, normal ring of characteristic $p$
with a perfect residue field and $\dim R=d$.
Let $M$ be a finitely generated $R$-module.
We show that
there exists $\beta(M) \in \mathbb R$ such that
$\lambda(M/I^{[q]}M) = e_{HK}(M) q^d + \beta(M) q^{d-1} + O(q^{d-2})$.
\end{abstract}

\maketitle

\section*{Introduction}
Throughout this paper $(R,\m)$ is a Noetherian local
${\mathbb Z}/p{\mathbb Z}$-algebra, $d=\dim R$,
and $I \subset R$ is
$\m$-primary.  We will let $n$ be a varying non-negative integer, and let
$q = p^n$. By  $I_n$ we
will denote the ideal generated by $x^{q}$,
$x\in I$.  If $M$ is a finite $R$-module, $M/I_nM$ has finite length;
we denote this length by $e_n(M,I)$, or more briefly by $e_n(M)$. We use
$\lambda(-)$ to denote
the length of an $R$-module. Let $f,g\colon {\mathbb N}\rightarrow {\mathbb
R}$ be functions from the nonnegative integers to the real numbers.
Recall that  $f(n)= O(g(n))$ if there exists a positive
constant $C$ such that $|f(n)|\leq Cg(n)$ for all $n\gg 0$, and
we write $f(n)=o(g(n))$ if $\lim_{n \rightarrow \infty} f(n)/g(n)=0$.

The basic question this paper studies is how
$e_n(M)$ depends on $n$.  The results of \cite{Mo1} show that
$e_n(M)=\alpha q^{d}+O(q^{d-1})$ for some real $\alpha$.  In section 1 we
strengthen this  proving:

\begin{T1}\label{Mtheorem1}
Let $(R,\m,k)$ be an excellent, local, normal ring of characteristic $p$
with a perfect residue field and $\dim R=d$.
Then $e_n(M)=\alpha q^{d}+\beta q^{d-1} +O(q^{d-2})$ for some $\alpha$ and
$\beta$ in $\mathbb R$.
\end{T1}

In the situation of Theorem 1\ref{Mtheorem1} it sometimes happens that
$\beta(M) =0$
provided that $M=R$ (or more generally when $M$ is torsion-free).
Our results establish that $\beta(M) =0$
whenever $M$ is torsion-free and the class group of $R$ is torsion.
In particular this holds when
$(R,\m,k)$ is a complete normal two-dimensional ring and $k$ is finite-see
Corollary 2.2.

One might hope that Theorem 1 could be generalized to prove that there
exists
a constant $\gamma$ such that
$e_n(M)=\alpha q^{d}+\beta q^{d-1} + \gamma q^{d-2} +O(q^{d-3})$ whenever
$R$ is non-singular in codimension two.  However, this cannot be true.
For example, if $R = {\mathbb Z}/5{\mathbb
Z}[x_1,x_2,x_3,x_4]/(x_1^4+\cdots+x_4^4)$,
then with $I = (x_1,...,x_4)$, $e_n(R) = \frac{168}{61}(5^{3n}) - \frac
{107}{61}(3^n)$ by
\cite{HaMo}. Note that $R$ is a 3-dimensional Gorenstein ring with isolated
singularity.
(Also see \cite{BuCh} for computations of the Hilbert-Kunz function
for plane cubics, as well as \cite{Mo2}-\cite{Mo4},
\cite{fakhruddintrivedi}, and \cite{Teixeira}
for other concrete computations of the Hilbert-Kunz function,
and \cite{WY1}-\cite{WY3}, \cite{BE} for work on minimal possible
values for the Hilbert-Kunz multiplicity.)
It may however be true that $e_n(M)$ is always
$\alpha q^d+\beta q^{d-1}+ (\mbox{periodic})(q^{d-2})+o(q^{d-2})$.
In particular we suspect that 
$e_n(M)= \alpha q^2+\beta q+(\mbox{eventually periodic})$ when
$d=2$.
Striking results in this
direction have recently been obtained by Brenner; see the
footnote in section 2. See also [Co] for the case $A=k[[x,y]]$, $I=(x,y)$,
$M$ arbitrary.
\section{ }
We make use of various facts about divisor classes in integrally
closed Noetherian domains.  Our reference is \cite{Bourbaki}, and we shall
need in particular Proposition 18 and Theorem 6 of chapter VII, section
4.

Let $R$ be an integrally closed Noetherian domain.  A Weil divisor on
$R$ is an element of the free abelian group on the height 1 primes of
$R$.  A principal Weil divisor is a divisor of the form $\sum_P
\ord_P(f)\cdot P$ with $f\neq 0$ in the field of fractions of $R$.
$C(R)$ is the quotient of the group of Weil divisors by the subgroup of
principal divisors.  Let $M$ be a finite $R$-module.  Then $M$ admits
a filtration with quotients (isomorphic to) $R/P_i$ where each $P_i$
is prime.  Consider the Weil divisor $-\sum P_i$, the sum extending
over those $P_i$ that are of height 1.  The image of this divisor in
$C(R)$ is independent of the choice of filtration, and is denoted by
$c(M)$.  The map $c$ is additive on exact sequences and $c(R)=0$.  If $P$ is a
height 1 prime of $R$ the exact sequence $0 \rightarrow P \rightarrow
R \rightarrow R/P \rightarrow 0$ shows that $c(P)=P$.  Suppose now that
we are in the situation of Theorem 1\ref{Mtheorem1} of the introduction.

\begin{lemma}\label{torsiontor1}
Let $(R,\m,k)$ be a  local ring
of characteristic $p$.
If $T$ is a finitely generated torsion $R$-module with $\dim T = \ell$,
then
$\lambda(\tor^{R}_1(R/I_n,T))=O(q^{\ell})$.
\end{lemma}

\begin{proof} Set $d = \dim R$.
Choose a system of parameters $(x_1, \ldots, x_d) \subseteq I$.
We induct on $\lambda(I/(x_1, \ldots, x_d))$.
If $\lambda(I/(x_1, \ldots, x_d))>0$, then there exists $J \subset I$
with $\lambda(I/J)=1$ so that we may write $I=(J,u)$ with $J: u = \m$.
For every $q = p^n$ there is an exact sequence
\[
0 \rightarrow
R/J_n\col u^q \rightarrow
R/J_n \rightarrow
R/I_n \rightarrow
0.
\]
Tensor with $T$ and look at the following portion of the long exact
sequence of Tors:
\[
\cdots \rightarrow
\tor^{R}_1(R/J_n,T) \rightarrow
\tor^{R}_1(R/I_n,T) \rightarrow
\tor^{R}_0(R/J_n\col u^q,T) \rightarrow \cdots.
\]
We have
$\lambda(\tor^{R}_1(R/J_n,T))=O(q^{\ell})$ by induction.
Also, since $J\col u=\m$, we have $\m_n\subseteq J_n\col u^q$ and
$\lambda(\tor^{R}_0(R/J_n\col u^q,T)) \leq
\lambda(\tor^{R}_0(R/\m_n,T))$.
But $\lambda(\tor^{R}_0(R/\m_n,T))$ is the Hilbert-Kunz function for
$T$, so $\lambda(\tor^{R}_0(R/\m_n,T))=O(q^{\ell})$.

We have reduced to the case where $\lambda(I/(x_1, \ldots, x_d))=0$.
We need a theorem which is implicitly in Roberts \cite{Roberts} and
explicitly given as Theorem 6.2 in  \cite{phantom} (see also
\cite[p278]{Seibert} for a theorem quickly giving an alternative proof
with a sharper result on the growth of the size of the Koszul groups):

\begin{nonumthm}
Let $(R,m)$ be a local ring of
characteristic $p$ and let $G_{\bullet}$ be a finite complex
$$
0\rightarrow G_{n}\rightarrow\dots\rightarrow G_{0}\rightarrow 0
$$
of length $n$ such that each $G_{i}$ is a finitely generated free module
and suppose that each $H_{i}(G_{\bullet})$ has finite length.  Suppose
that $M$ is a finitely generated $R$-module.  Let $ d = \dim  M.$  Then
there is a constant $C > 0$ such that $\ell(H_{n-t}(M\otimes_{R}
F^{e}(G_{\bullet}))\leq Cq^{\min(d,t)}$ for all $t\geq 0$ and all $e\geq
0$, where $q =
p^{e}$.
\end{nonumthm}

Consider $K_{\bullet} ((\underline x);R)$, the Koszul complex on $(x_1,
\ldots, x_d)$.  Let $H_{\bullet}((\underline x);R)$ denote the
homology of the Koszul complex. We apply the above theorem to conclude
that there exists a
constant $C>0$ such that
$\lambda(H_{d-t}(T\otimes F^n(K_{\bullet})))\leq Cq^{\min(\ell,t)}$ for
all $t$ and for all $n$.
Hence $\lambda(H_i(T\otimes F^n(K_{\bullet})))\leq O(q^{\ell})$ for all $i$.
In general, $H_1(T\otimes F^n(K_{\bullet}))$ maps onto
$\tor^{R}_1(T,R/I_n)$, which gives the stated result.
(To see this, note that $F^n(K_{\bullet})$ is exactly the Koszul complex
on the generators of $I$ raised to the $q = p^n$ power, so that both
of the complexes $F^n(K_{\bullet})$ and the minimal free resolution of
$R/I_n$
begin with the same two free modules.)
\end{proof}

\begin{lemma}\label{Mlemma1.1}
Let $(R,\m,k)$ be a  local, normal ring
of characteristic $p$.
Let $J\neq (0)$ be an ideal of $R$.  If $c(J)=0$ then
$e_n(J)=e_n(R)+O(q^{d-2})$.
\end{lemma}

\begin{proof}

The primary decomposition theorem shows that $J= (\cap P_i^{(n_i)})
\cap N$ where the $P_i^{(n_i)}$ are symbolic powers of finitely many
height one
primes $P_i$, and $\dim(R/N)\leq d-2$.  Then $c(J)$ is the divisor class
represented by $\sum n_i P_i$.  So $\sum n_i P_i$ is the divisor of
some $f\neq 0$ in $R$, and replacing $J$ by $f^{-1}J$ we may assume
that $T=R/J$ is of dimension at most $d-2$.
By Lemma \ref{torsiontor1},
$\lambda(\tor_1^R(T, R/I_n)) =  O(q^{d-2})$ since $\dim T \leq d-2$.
The exact
sequence
$\tor^{R}_1(R/I_n,T)\rightarrow J/I_n J \rightarrow R/I_n \rightarrow
T/I_nT$ now gives the result, since $\lambda(T/I_nT) = O(q^{d-2})$ by
\cite{Mo1}.
\end{proof}

\begin{lemma}\label{Mlemma1.2}
Let $(R,\m,k)$ be a local, normal ring
of characteristic $p$.
Suppose $M$ is torsion-free of rank $r$ and $c(M)=0$.  Then
$e_n(M)\leq re_n(R)+O(q^{d-2})$.
\end{lemma}

\begin{proof}
By Theorem 6, section 4 of \cite{Bourbaki} there is an exact sequence
$0 \rightarrow R^{r-1} \rightarrow M \rightarrow J \rightarrow 0$ \
with
$J$ an ideal.  Then $c(J)=0$, and we use the exact sequence
$(R/I_n)^{r-1}\rightarrow M/I_nM \rightarrow J/I_nJ$
together with
Lemma \ref{Mlemma1.1}.
\end{proof}

\begin{theorem}\label{Mtheorem1.3}
Let $(R,\m,k)$ be a local, normal ring
of characteristic $p$.
Suppose $M$ is torsion-free of rank $r$ and $c(M)=0$.  Then $e_n(M)
=re_n(R) +O(q^{d-2})$.
\end{theorem}

\begin{proof}
There is an exact sequence
$0 \rightarrow K \rightarrow R^{r+s} \rightarrow M \rightarrow 0$
for some $K$ and $s\geq 0$.  Then $K$ is torsion-free of rank $s$ and
$c(K)=0$.  By Lemma \ref{Mlemma1.2}, $e_n(K)\leq se_n(R)+O(q^{d-2})$.
Evidently $e_n(K)+e_n(M)\geq (r+s)e_n(R)$.  So $e_n(M) \geq
re_n(R)+O(q^{d-2})$;  Lemma \ref{Mlemma1.2} provides the opposite
inequality.
\end{proof}

\begin{lemma}\label{Mlemma1.4}
Let $(R,\m,k)$ be a local, normal ring
of characteristic $p$.
If $M$ is torsion-free, $\lambda(\tor^{R}_1(R/I_n,M)) = O(q^{d-2})$.
\end{lemma}

\begin{proof}
There is an exact sequence
$0 \rightarrow K \rightarrow G \rightarrow M \rightarrow 0$
with $G$ free.  Then $c(K\oplus M)=c(G)=0$, and Theorem
\ref{Mtheorem1.3} applied to $K\oplus M$ shows that
$e_n(G)=e_n(K)+e_n(M)+O(q^{d-2})$.  Now use the exact sequence
$0 \rightarrow \tor^{R}_1(R/I_n,M) \rightarrow K/I_nK \rightarrow G/I_nG
\rightarrow M/I_nM \rightarrow 0$.
\end{proof}

\begin{lemma}\label{Mlemma1.5}
Let $(R,\m,k)$ be a local, normal ring
of characteristic $p$.
Suppose $M$ and $N$ are torsion-free of the same rank and
$c(M)=c(N)$.  Then $e_n(M)=e_n(N)+O(q^{d-2})$.
\end{lemma}

\begin{proof}
Replacing $M$ and $N$ by $M\oplus J$ and $N\oplus J$ for some ideal
$J$ we may assume that $c(M)=c(N)=0$.  Now apply Theorem
\ref{Mtheorem1.3} to $M$ and to $N$.
\end{proof}

\begin{definition}\label{Mdef1.6}
Let $(R,\m,k)$ be a local, normal ring
of characteristic $p$.
If $M$ is torsion-free of rank $r$, $\delta_n(M)=e_n(M)-re_n(R)$.
\end{definition}

\begin{rmks}
$\delta_n(R)=0$ and $\delta_n(M\oplus N)=\delta_n(M)+\delta_n(N)$.  If
$c(M)=c(N)$, Lemma \ref{Mlemma1.5} tells us that
$\delta_n(M)=\delta_n(N)+O(q^{d-2})$.
\end{rmks}

To make further progress we shall use the $p$th power map
$F\colon R\rightarrow R$, assuming $R$ is complete with perfect residue field. In this case $F$ is finite of degree $p^d$.
Given a finite map $R \rightarrow R'$ between integrally closed
Noetherian domains, we obtain induced norm maps from Weil divisors on
$R'$ to Weil divisors on $R$ and from $C(R')$ to $C(R)$.  For
$F\colon R \rightarrow R$
we claim that these norm maps are just multiplication by
$p^{d-1}$.  For if $P$ is a height 1 prime of $R$, the only prime lying over
$P$ is $P$ itself, and the ramification degree is evidently $p$.  So
the residue class field degree is $p^{d-1}$ by the discussion of Section
4.8, Chapter VII, page 535 in
\cite{Bourbaki}, and then the norm of $P$ is
$p^{d-1}\cdot P$ by the same discussion.

If $M$ is a finitely generated $R$-module of rank $r$ let $^1M$ be
$M$ as additive group, but with $R$ acting through the $p$th power map
$F\colon R \rightarrow R$.  Then $^1M$ is evidently finite of rank
$p^dr$, and $e_n(^1M)=e_{n+1}(M)$ for all $n$.

\begin{theorem}\label{Mtheorem1.7}
Let $(R,\m,k)$ be an excellent, local, normal ring
of characteristic $p$
with
a perfect residue field.
Let $M$ be torsion-free of rank $r$.
Then $\delta_{n+1}(M)=p^{d-1}\delta_n(M)+O(q^{d-2})$.
\end{theorem}

\begin{proof} We may complete $R$ without changing the hypotheses or
conclusions,
and henceforth we assume that $R$ is complete.
Since the norm map $C(R)\rightarrow C(R)$ induced by $F\colon R \rightarrow
R$ is multiplication by $p^{d-1}$, Proposition 18, section 4.8, Chapter
VII of \cite{Bourbaki}
tells us that
$c(^1M)=p^{d-1}c(M)+rc(^1R)$.
The remarks after Definition \ref{Mdef1.6} then show that
$\delta_n(^1M)-r\delta_n(^1R)=p^{d-1}\delta_n(M)+O(q^{d-2})$.
But
$\delta_n(^1M)-r\delta_n(^1R)=
(e_{n+1}(M)-p^dre_n(R))-r(e_{n+1}(R)-p^de_n(R))=
e_{n+1}(M)-re_{n+1}(R)=
\delta_{n+1}(M)$,
giving the theorem.
\end{proof}

\begin{theorem}\label{Mtheorem1.8}
Let $(R,\m,k)$ be an excellent, local, normal ring
of characteristic $p$
with
a perfect residue field. Let $M$ be a torsion-free finite $R$-module.
There is a real constant $\tau(M)$ such that
$\delta_n(M)=\tau(M)q^{d-1}+O(q^{d-2})$.
\end{theorem}

\begin{proof}
Let $v_n=\delta_n(M)/q^{d-1}$.  By Theorem \ref{Mtheorem1.7},
$v_{n+1}-v_n=O(1/q)$.  So $v_n \rightarrow$ some $\tau$, and
$v_n =\tau + O(1/q)$. The result follows.
\end{proof}

\begin{cor}\label{Mcor1.9}
Let $(R,\m,k)$ be an excellent, local, normal ring
of characteristic $p$
with
a perfect residue field.
There is a homomorphism $\tau\colon C(R)\rightarrow \mathbb R, +$ with the
following property.  If $M$ is torsion-free of rank $r$ then
$e_n(M)=re_n(R)+\tau q^{d-1} +O(q^{d-2})$ with $\tau = \tau(c(M))$.
\end{cor}

\begin{proof}
If $c\in C(R)$ choose $M$ torsion-free with $c(M)=c$.  Then
$\delta_n(M)=\tau q^{d-1} +O(q^{d-2})$ for some real $\tau$.  The remarks
after Definition \ref{Mdef1.6} tell us that $\tau$ is independent of
the choice of $M$ and that $c \rightarrow \tau$ is a homomorphism.
\end{proof}

We remark that it is immediate from this corollary that $\tau$ is the
zero map whenever the class group of $R$ is torsion.

\begin{theorem}\label{thmforR}
Let $(R,\m,k)$ be an excellent, local, normal ring
of characteristic $p$
with
a perfect residue field.  Let $\dim R=d$.
Then there exists $\beta(R) \in \mathbb R$ such that
$e_n(R) = e_{HK}(I;R)q^d + \beta(R) q^{d-1} + O(q^{d-2})$.
Furthermore, $\beta(R)(p^{d-1}-p^d) =
\tau(^1R)$.
\end{theorem}

\begin{proof}
Taking $M={}^1R$ in Theorem \ref{Mtheorem1.8} we find that     $e_{n+1}(R)-p^de_n(R)=\tau q^{d-1} +O(q^{d-2})$ where     $\tau=\tau(^1R)$.  Set $u_n=e_n(R)-\beta q^{d-1}$ where
$(p^{d-1}-p^{d})\beta = \tau$.  Then
$u_{n+1}-p^du_n=e_{n+1}(R)-p^de_n(R)-\tau q^{d-1} =O(q^{d-2})$, and arguing
as in the proof of Theorem \ref{Mtheorem1.8} we find that $u_n=\alpha
q^{d}+O(q^{d-2})$.

In other words,
$e_n(R) = \alpha(R)q^d + \beta(R) q^{d-1} + O(q^{d-2})$ where
$\beta(R)=\tau({}^1R)/(p^{d-1}-p^d)$. Clearly $\alpha(R) =
e_{HK}(I;R)$ is forced.
\end{proof}

\begin{theorem}\label{mainthm}
Let $(R,\m,k)$ be an excellent, local, normal ring of characteristic $p$
with a perfect residue field and $\dim R=d$.
Let $M$ be finitely generated $R$-module.
Then there exists $\beta(M) \in \mathbb R$ such that
$e_n(M) = e_{HK}(I;M)q^d + \beta(M) q^{d-1} + O(q^{d-2})$.
\end{theorem}

\begin{proof}
We again complete $R$ and assume it is complete.
Suppose
first that $M$ is torsion-free.  Then the result follows from Theorems
\ref{Mtheorem1.8} and \ref{thmforR}.  In general there is an exact
sequence
$0\rightarrow T \rightarrow M \rightarrow M' \rightarrow 0$
with $T$ torsion and $M'$ torsion free.  The exact sequence
$\tor^{R}_1(R/I_n,M')\rightarrow T/I_nT\rightarrow M/I_nM\rightarrow
M'/I_nM' \rightarrow 0$ combined with Lemma \ref{Mlemma1.4} shows
that
$e_n(M)=e_n(M')+e_n(T)+O(q^{d-2})$.  Since $T$ has dimension $\leq
d-1$, \cite{Mo1} shows that $e_n(T)=cq^{d-1}+O(q^{d-2})$ for some $c\geq
0$, and the
result for $M'$ yields the result for $M$.
\end{proof}

A corollary of the above results gives us similar growth conditions on
certain Tor modules.

\begin{cor}\label{cortorsion}
Let $(R,\m,k)$ be an excellent, local, normal ring
of characteristic $p$ with perfect residue field and
with $\dim R=d$.
Let $T$ be a torsion $R$-module.
Then there exists $\gamma(T) \in \mathbb R$ such that
$\lambda(\tor_1^R(T, R/I_n))  = \gamma(T) q^{d-1} + O(q^{d-2})$.
\end{cor}

\begin{proof} We may complete $R$ and henceforth assume $R$ is complete.
Consider an exact sequence,
\[
0\rightarrow M\rightarrow R^s\rightarrow T\rightarrow 0
\]
where $M$ is torsion free. The long exact sequence on Tor after
tensoring with
$R/I_n$ shows that
$$\lambda(\tor_1^R(T, R/I_n)) = e_n(M) + e_n(T) - se_n(R).$$
By \cite{Mo1}, $e_n(T) = cq^{d-1}+O(q^{d-2})$ for some $c\geq 0$, while
Theorem \ref{Mtheorem1.8} shows that $e_n(M) - se_n(R) = \tau(M)q^{d-1} +
O(q^{d-2})$.
The corollary follows.
\end{proof}
\bigskip
\section{The map $\tau$} 
\medskip

Given $(R,\m)$ and $I$ as in the last section, it might  seem plausible that the
map $\tau\colon C(R)\rightarrow \mathbb R$ is always the zero map, so that $e_n(M)=cq^d+O(q^{d-2})$
for torsion-free $M$. This of course is true when $R$ is a UFD (or more generally
when $C(R)$ is torsion), and so holds for $k[[x_1,...,x_r]]/(F)$ when $F$ is a
smooth form and $r$ is at least 5.
However,  K.-i. Watanabe has found counterexamples when $d$ is $3$ or more. Here's
a very nice example of his. Let $S$ be the $2\times 3$ matrix whose entries are independent
variables $x_1,...,x_6$, and let $R$ be the quotient of $k[[x_1,...,x_6]]$ by the 
ideal generated by the $2\times 2$ minors of $S$. If $I = \m$, then $e_n(R)=
(13q^4-2q^3-q^2-2q)/8$. Hence $\beta(R)=-1/4$, and $\tau(^1 R)=(p^4-p^3)/4$.
K. Kurano also commented on this point to us. He reports that
$\beta(R) = 0$ if the canonical
class $c(K_R)$ is torsion in the divisor class group $C(R)$. (Here $R$ is local, excellent
and normal.) His proof uses the singular Riemann-Roch theorem, and 
furthermore shows that if $R$ is Cohen-Macaulay and $I$ is a maximal primary
ideal of finite projective dimension,   then $e_n(R)$ is a polynomial
in $q$ with rational coefficients.\footnote{A recent preprint of
  H. Brenner \cite{Br1} shows that the
  Hilbert-Kunz multiplicity of the ring
is rational in the two-dimensional graded case.  This result was
obtained independently by V. Trivedi \cite{Trivedi}.
In  another more recent preprint \cite{Br2} Brenner proves
that $e_n(R) = \alpha q^2 + \mbox{ an eventually
  periodic function of }n$
in the two-dimensional graded case over the algebraic closure of a
finite field.}

When $d=2$ there are more general results. In particular
the following Lemma seems to be known to experts, and we
thank M. Artin and J. Lipman for pointing out relevant references and facts.

\begin{lemma}\label{torsionclassgroup}
Suppose that  $(R,\m,k)$  is a complete local normal two-dimensional ring, and $k$ is the algebraic closure of the
field with $p$ elements.  Then $C(R)$ is a torsion group.
\end{lemma}

\begin{proof} The proof depends on the numerical theory of exceptional
divisors (treated in full generality by Lipman), and arguments of Artin. An
exposition is given by H. G\"ohner in \cite{Go},
section 4, pages 423-426, which is independent of the rest of G\"ohner's paper.
Note in particular the first part of Theorem 4.4
and corollary 4.5 in this paper. The hypothesis that there is a desingularization
$f\colon X\rightarrow \text{Spec}(R)$, made at the beginning of section 4, is satisfied in this case, see \cite{Li2}. \end{proof}
\medskip
\begin{cor}\label{tauiszero} Suppose that  $(R,\m,k)$  is a complete local normal two-dimensional ring, and $k$ is finite.
Then $\tau$ is the zero map.
\end{cor}

\begin{remark} For general algebraically closed $k$ there is an analog of Lemma \ref{torsionclassgroup}.
We adopt the notation of \cite{Go}. By (*) on page 425 there is an exact
sequence $$0\rightarrow \text{Pic}^0(X)\rightarrow C(R)\rightarrow H\rightarrow 0$$ with $H$ finite; see page 425 for
the definition of $\text{Pic}^0(X)$. To prove Lemma \ref{torsionclassgroup}, G\"ohner uses Artin's
result that there is a filtration of $\text{Pic}^0(X)$ with each quotient
isomorphic to either the additive group of  $k$,  $k^*$, or the group of $k$-valued points of
the Jacobian variety of an irreducible component of the reduced special
fibre of $f$. Somewhat more is true. There is
a connected algebraic group $G$ defined over $k$, built out of copies of the
additive group, the multiplicative group and the above Jacobians, such
that $\text{Pic}^0(X)$ identifies with $G_k$. For more information
concerning this topic, see \cite{Li1}, in particular Theorem 7.5.
\end{remark}

\begin{remark}    We believe that Corollary \ref{tauiszero}  holds even when $k$ is
infinite. Here's an intuitive argument. Suppose that $P$ and $Q$
are in some sense ``generic points'' of $G_k = \text{Pic}^0(X)$. Because the
definition of $\tau$ is purely algebraic, $\tau(P)= \tau(Q)$. Since the various
$P-Q$ with $P$ and $Q$ generic generate $G_k$, at least when $k$ is large enough, $\tau$
vanishes on the subgroup $\text{Pic}^0(X)$ of $C(R)$ of finite index.
\end{remark}

\begin{remark} The third author has made the idea of the above remark into a simple proof when
$R$ is the homogeneous coordinate ring of a smooth projective curve,
localized at the homogeneous maximal ideal. In particular when $R=
k[[x_1,x_2,x_3]]/(F)$, $F$ a smooth form, $\tau$ is the zero map. As we've noted
this is also true for 5 or more variables--the 4 variable case remains
open.

\end{remark}


\begin{thebibliography}{HaMo}

\bibitem[BE]{BE}
M.~Blickle and F.~Enescu.
\newblock On rings with small Hilbert-Kunz functions.
\newblock {\em Proc. Amer. Math. Soc.},  132(9): 2505--2509, 2004.

\bibitem[Bo]{Bourbaki}
N.~Bourbaki.
\newblock {\em Elements of mathematics. {C}ommutative algebra}.
\newblock Hermann, Paris, 1972.
\newblock Translated from the French.

\bibitem[Br1]{Br1}
H.~Brenner.
\newblock The rationality of the Hilbert-Kunz multiplicity in graded dimension two.
\newblock Preprint, (2004), \ttfamily arXiv:math.AC/0402180\normalfont. 

\bibitem[Br2]{Br2}
H.~Brenner.
\newblock The Hilbert-Kunz function in graded dimension two.
\newblock Preprint, (2004), \ttfamily arXiv:math.AC/0405202\normalfont.

\bibitem[BuCh]{BuCh}
R.-O. Buchweitz and Q. Chen.
\newblock Hilbert-Kunz functions of cubic curves and surfaces.
\newblock {\em J. Algebra}, 197(1): 246--267, 1997.

\bibitem[Co]{Contessa}
M.~Contessa.
\newblock On the {H}ilbert-{K}unz function and {K}oszul homology.
\newblock {\em J. Algebra}, 175(3): 757--766, 1995.

\bibitem[FaTr]{fakhruddintrivedi}
N.~Fakhruddin and V.~Trivedi.
\newblock Hilbert-{K}unz functions and multiplicities for full flag varieties
  and elliptic curves.
\newblock {\em J. Pure Appl. Algebra}, 181(1): 23--52, 2003.

\bibitem[G\"o]{Go}
H.~G\"ohner.
\newblock Semifactoriality and Muhly's condition (N) in two dimensional
local rings.
\newblock {\em J. Algebra}, 34: 403--429, 1975.

\bibitem[HH]{phantom}
M.~Hochster and C.~Huneke.
\newblock Phantom homology.
\newblock {\em Mem. Amer. Math. Soc.}, 103(490):vi+91, 1993.

\bibitem[HaMo]{HaMo}
C.~Han and P.~Monsky.
\newblock Some surprising {H}ilbert-{K}unz functions.
\newblock {\em Math. Z.}, 214(1):119--135, 1993.

\bibitem[Li1]{Li1}
J.~Lipman.
\newblock The Picard group of a scheme over an Artin ring.
\newblock {\em Inst. Hautes \'Etudes Sci. Publ. Math.}, 46: 15--86, 1976.

\bibitem[Li2]{Li2}
\bysame
\newblock Desingularization of two-dimensional schemes.
\newblock {\em Annals of Math.}, 107: 151--207, 1978.

\bibitem[Mo1]{Mo1}
P.~Monsky.
\newblock The {H}ilbert-{K}unz function.
\newblock {\em Math. Ann.}, 263(1): 43--49, 1983.

\bibitem[Mo2]{Mo2}
\bysame
\newblock The Hilbert-Kunz function of a characteristic $2$ cubic.
\newblock {\em J. Algebra}, 197(1): 268--277, 1997.

\bibitem[Mo3]{Mo3}
\bysame
\newblock Hilbert-Kunz functions in a family: point-$S\sb 4$ quartics.
\newblock {\em J. Algebra}, 208(1): 343--358, 1998.

\bibitem[Mo4]{Mo4}
\bysame
\newblock Hilbert-Kunz functions in a family: line-$S\sb 4$ quartics.
\newblock {\em J. Algebra}, 208(1): 359--371, 1998.

\bibitem[Ro]{Roberts}
P.~Roberts.
\newblock Le th\'eor\`eme d'intersection.
\newblock {\em C. R. Acad. Sci. Paris S\'er. I Math.}, 304(7):177--180,
1987.

\bibitem[Se]{Seibert}
G.~ Seibert.
\newblock Complexes with homology of finite length and {F}robenius functors.
\newblock {\em J. Algebra}, 125(2):278--287, 1989.

\bibitem[Te]{Teixeira}
P.~Teixeira
\newblock $p$-fractals and Hilbert-Kunz series.
\newblock Ph.D. thesis, Brandeis University, 2002.

\bibitem[Tr]{Trivedi}
V.~Trivedi.
\newblock Semistability and {H}ilbert-{K}unz multiplicity for curves.
\newblock Preprint, (2004), \ttfamily arXiv:math.AC/0402245\normalfont.

\bibitem[WY1]{WY1}
K.-i.~Watanabe and K.-i.~Yoshida.
\newblock Hilbert-Kunz multiplicity and an inequality between multiplicity and colength.
\newblock {\em J. Algebra}, 230(1): 295--317, 2000.

\bibitem[WY2]{WY2}
\bysame
\newblock Hilbert-Kunz multiplicity of two-dimensional local rings.
\newblock {\em Nagoya Math. J.}, 162: 87--110, 2001.

\bibitem[WY3]{WY3}
\bysame
\newblock Hilbert-Kunz multiplicity of three-dimensional local rings.
\newblock Preprint, (2003), \ttfamily arXiv:math.AC/0307294\normalfont.

\end{thebibliography}

\end{document}